\documentclass[12pt,a4paper]{article}

\usepackage{latexsym}
\usepackage{amsfonts}
\usepackage{amsmath}
\usepackage{amssymb}
\usepackage{color}
\usepackage{theorem}
\usepackage[dvips,
colorlinks=true,
linkcolor=webred,
filecolor=webred,
citecolor=webred,
pdftitle={Groups of prime-power order with a small second derived quotient},
pdfauthor={Csaba Schneider},
bookmarksopen=true]{hyperref}

\definecolor{webred}{cmyk}{0,1,.7,.3}
\definecolor{webgreen}{cmyk}{1,0,1,.4}
\definecolor{webblue}{cmyk}{1,1,0,.5}

\theorembodyfont{\sl}

\newtheorem{theorem}{Theorem}[section]

\newtheorem{lemma}[theorem]{Lemma}
\newtheorem{corollary}[theorem]{Corollary}

\theorembodyfont{\rm}

\newtheorem{example}[theorem]{Example}

\newenvironment{proof}{{\sc Proof.}}{\hfill$\Box$

\bigskip

}

\newcommand{\der}[2]{#1^{(#2)}}
\newcommand{\lcs}[2]{\gamma_{#2}(#1)}

\newcommand{\gen}[1]{\bigl<#1\bigr>}

\newcommand{\cent}[2]{{\sf C}_{#1}\left(#2\right)}

\newcommand{\replace}{\leadsto}

\makeatletter
\addto@hook\every@math@size{\dch@script@adjust}
\def\dch@script@adjust{\@ifundefined{dch@sizet\f@size}%
  {\expandafter\dch@set@script\csname dch@sizet\f@size\endcsname}%
  {\csname dch@sizet\f@size\endcsname}}
\def\dch@set@script#1{\begingroup
  \frozen@everymath{}
  \let#1\@empty \let\dch@do@one\relax
  \dch@set@one\scriptscriptstyle \scriptscriptfont#1\ssf@size
  \dch@set@one\scriptstyle\scriptfont#1\sf@size
  \dch@set@one\textstyle\textfont#1\f@size
  \endgroup #1} %
\def\dch@set@one#1#2#3#4{%
  \@ifundefined{dch@size#4}%
   {\expandafter\xdef\csname dch@size#4\endcsname{%
      \fontdimen13\the#2\tw@=\the\fontdimen13#2\tw@
      \fontdimen14\the#2\tw@=\the\fontdimen14#2\tw@
      \fontdimen15\the#2\tw@=\the\fontdimen15#2\tw@
      \fontdimen16\the#2\tw@=\the\fontdimen16#2\tw@
      \fontdimen17\the#2\tw@=\the\fontdimen17#2\tw@}%
  }{\csname dch@size#4\endcsname}%
  \setbox\z@\hbox{$#1H_2$}\@tempdima\dp\z@
  \setbox\z@\hbox{$#1H_2^{+\vrule \@height 1em}$}%
   \ifdim\@tempdima<\dp\z@
    \advance\@tempdima\dp\z@ \divide\@tempdima\tw@
    \@tempdimb\dp\z@ \advance\@tempdimb-\@tempdima
    \advance\@tempdimb\ht\z@ \advance\@tempdimb-1em
    \xdef#3{#3\dch@do@one#2{\the\@tempdimb}{\the\@tempdima}}%
  \fi}
\def\dch@do@one#1#2#3{\fontdimen13#1\tw@#2\relax
  \fontdimen14#1\tw@\fontdimen13#1\tw@ \fontdimen15#1\tw@\fontdimen13#1\tw@
  \fontdimen\sixt@@n#1\tw@#3\fontdimen17#1\tw@\fontdimen\sixt@@n#1\tw@}%
\makeatother

\begin{document}

\title{Groups of prime-power order with a small second derived quotient}
\author{Csaba Schneider\\
School of Mathematics and Statistics\\
The University of Western Australia\\
35~Stirling Highway
6009 Crawley\\
Western Australia\\
{\tt www.maths.uwa.edu.au/$\sim$csaba}\\
{\tt csaba@maths.uwa.edu.au}
\date{9 January 2003}}

{\small \maketitle}

\abstract 
For odd primes we prove some structure
theorems for finite $p$-groups $G$, such that $G''\neq 1$ and
$|G'/G''|=p^3$.
Building on results of Blackburn and Hall, it is shown that $\lcs G3$ is a maximal subgroup of
$G'$, the group $G$ has a central decomposition into two simpler
subgroups, and, moreover, 
$G'$ has  one of two isomorphism types.
\endabstract


\section{Introduction}\label{intro}

It is well-known that in a finite $p$-group $G$ the condition
$G''\neq 1$ implies that $|G'/G''|\geqslant p^3$;  see for example Huppert~\cite{Huppert}~III.7.10.  In this article we
prove a number of results about groups in which equality holds; that is, we assume that $G''\neq 1$ and
$|G'/G''|=p^3$. Such groups have already been investigated by, among
others, N.~Blackburn and P.~Hall. Blackburn~\cite{Blackburn84} proved
that the condition $|G'/G''|=p^3$ implies that $G''$ is abelian
generated by two elements and it is nearly homocyclic. 
In the same article he also published a result, which he attributed
to Hall, that for odd primes the same condition implies that
$|G''|\leqslant p$. Here we mostly consider $p$-groups for odd
$p$, and our main results are concerned with such groups.

Let $G$ be a finite $p$-group and $\lcs Gi$ the $i$-th term of the lower central series, so
that $\lcs G1=G$, $\lcs G2=G'$, etc. If $G''\neq 1$ then we have the following chain of normal subgroups:
\begin{equation}\label{normsgps}
G> G'=\lcs{G}{2}>\lcs{G}{3}>\lcs{G}{4}\geqslant G''>
1.
\end{equation}
If, in addition, we assume that $|G'/G''|=p^3$, then it easily follows 
that the order of $G'/\lcs{G}{3}$ is at most $p^2$. The result of this
simple argument is improved by the following theorem.

\begin{theorem}\label{3pgr}
Let $p\geqslant 3$ and $G$ be a finite $p$-group, such that
$|G'/G''|=p^3$ and $G''\neq 1$. 
Then $|G'/\lcs{G}{3}|=p$ and
$G''=\lcs{G}{5}$.
\end{theorem}

The proof of this result is given in Section~\ref{sec2}. 
Our second theorem, whose proof is in Section~\ref{sec3},
is that $G$ can be written as a central product of two simpler subgroups.

\begin{theorem}\label{decomp2}
Let $p\geqslant 3$ and $G$ be a finite $p$-group, such that
$|G'/G''|=p^3$ and $G''\neq 1$. Then $G$ can be factorised as 
$G=HU$, where
\begin{enumerate}
\item[(i)] $H$ is a normal subgroup of $G$ generated by at most $5$
generators;
\item[(ii)] $\lcs Hi=\lcs Gi$ for all $i\geqslant 2$;
\item[(iii)] $U$ is a normal subgroup of $G$, such that $U'\leqslant\lcs{G}{5}$;
\item[(iv)] $H$ and $U$ centralise each other.
\end{enumerate}
\end{theorem}

An example is given after the proof of this theorem
to show that the number ``5'' is, in general, best possible, and 
that there are, in some cases, other central decompositions of $G$ 
in which the subgroups can have different isomorphism types.

Our proofs are based on commutator calculus.
To simplify notation, we write long commutators according to the left-normed convention; for example $[a,b,c]=[[a,b],c]$. 
We use the well-known commutator identities that can be found in most
group theory textbooks (see for instance
Huppert~\cite{Huppert}~III.1.2-III.1.3). 
In addition to these, we need the collection formula,
which is proved as Lemma~VIII.1.1 by Huppert and
Blackburn~\cite{HBII}.
We mainly use this result in the simplest case when it can
be stated as
$$
[x^p,y]\equiv[x,y]^p\bmod (N')^p\lcs{N}{p}\quad\mbox{where}\quad
N=\gen{x,[x,y]}.
$$
The Hall-Witt identity will occur in a lesser known form which can be
found in Magnus, Karrass \& Solitar~\cite{Magnusetal} on page~290:
$$
[x,y,z^x][z,x,y^z][y,z,x^y]=[x,y,z[z,x]][z,x,y[y,z]][y,z,x[x,y]]=1.
$$

We often manipulate generating sets of groups. In order to avoid
cumbersome repetitions, we introduce a piece of notation. Let $G$ be a
group, $g$ a symbol referring to  a group element, and $x$ an element
in $G$. 
After the  occurrence of the expression $x\replace g$, the name $g$ will
refer to the element $x$. For example, let $G$ be the cyclic group of
order two and let $g$ denote its non-identity element. If we perform the
replacement $g^2\replace g$, then the symbol $g$ will refer to the
identity element of $G$. 

One can naturally ask whether it
is possible for a fixed prime to give a classification of groups
which satisfy the conditions of the previous two   theorems. 
It is conceivable that
Blackburn's~\cite{Blackburn54} description of groups of maximal class with order~$p^6$
and degree of commutativity~0 is  a good starting point. However,
increasing the number of generators and allowing the abelian factor to
have exponent higher than $p$ led to complications which could not
be resolved within the research presented here.

Our results can also be viewed in a wider context. It was first shown
by Hall~\cite{Hall34} (Theorem~2.57) that the conditions $i\geqslant 1$ and $\der
G{i+1}\neq 1$ imply that $|\der Gi/\der G{i+1}|\geqslant p^{2^i+1}$,
and
$|G|\geqslant p^{2^{i+1}+i+1}$ (see also Huppert~\cite{Huppert}
III.7.10 and III.7.11). The lower bound for the order of $G$
has recently been improved by Mann~\cite{mann} and the author~\cite{csthesis}. Both of these
improvements are, however, minor, and the order of the smallest
$p$-group $G$ such that $\der G{i+1}\neq 1$ is still unknown; the smallest known
examples were constructed by Evans-Riley, Newman, and the
author~\cite{se-rmfncs}. If $p\geqslant 3$ then we also do not know
how sharp Hall's lower bound is for $|\der Gi/\der G{i+1}|$.
As the example of the Sylow 2-subgroup of the symmetric group with degree
$2^{i+2}$ shows, this result is best possible for $p=2$; it is not known otherwise. Our research was originally motivated by
these questions, and it is hoped that a more detailed understanding of
groups with a small second derived quotient will give us a hint of the 
solution to some of the above problems. Some partial results can be
found in the author's PhD thesis~\cite{csthesis}.

Our results are inspired by Lie algebra calculations, and  it is possible
to prove some of them using the Lie ring method.
In fact, Theorem~\ref{3pgr} can be proved by
first verifying the corresponding result for Lie algebras and then using
the  Lie ring associated with the lower central series. This approach would
lead to some interesting new results for Lie algebras, which are
beyond the scope of the present article.

The paper is structured as follows. In Section~\ref{sec1} we prove a
lemma which is a generalisation of Blackburn's
Theorem~1.3~\cite{Blackburn54}. A consequence of this result is that
we can often restrict our interest to groups which are generated by
two or three elements.
In Sections~\ref{sec2} and~\ref{sec3}  we prove Theorems~\ref{3pgr}
and~\ref{decomp2}, respectively. 
In Section~\ref{sec4} we characterise the commutator subgroup of $G$, and
show that it has one of two isomorphism types. 

\section{A general lemma and some consequences}\label{sec1} 

We have seen in the introduction that in a finite $p$-group $G$, the
conditions $|G'/G''|=p^3$ and $G''\neq 1$ imply that $G'/\lcs G3$ has
order at most $p^2$. The aim of this section
is to show that $G$ has a subgroup $H$ with a small generating set,
such that, apart from the first term, the lower central series of $H$
coincides with the lower central series of $G$. This result
generalises  Blackburn's Theorem~1.3~\cite{Blackburn54} and
Slattery's Lemma~2.1~\cite{slat}. 

\begin{lemma}\label{gh}
Let $G$ be a nilpotent group and $H$ a subgroup of $G$, such that
$G'=H'\lcs G3$. Then $\lcs Gi=\lcs Hi$ for all $i\geqslant
2$. Moreover, $H$ is a normal subgroup of $G$.
\end{lemma}
\begin{proof}
First we prove by induction on $i$ that $\lcs Gi=\lcs Hi\lcs
G{i+1}$ for all $i\geqslant 2$. By the conditions of the lemma, this is true for
$i=2$. Suppose that our claim holds for some $i-1\geqslant 2$, and let us
show that it holds for $i$ as well. As it is obvious that $\lcs Hi\lcs
G{i+1}\leqslant \lcs Gi$, we only have to prove that $\lcs
Gi\leqslant\lcs Hi\lcs G{i+1}$. Using the induction hypothesis and
III.1.10(a) of Huppert~\cite{Huppert}, we compute
\begin{multline*}
\lcs G{i}=[\lcs G{i-1},G]=
[\lcs H{i-1}\lcs Gi,G]\\
={[\lcs
H{i-1},G]}[\lcs
Gi,G]=[\lcs H{i-1},G]\lcs G{i+1}.
\end{multline*}
Therefore it is enough to prove that $[\lcs H{i-1},G]\leqslant \lcs
H{i}\lcs G{i+1}$. First we note that 
$\lcs Gi\geqslant \lcs Hi\lcs G{i+1}\geqslant \lcs G{i+1}$,
and hence $\lcs Hi\lcs G{i+1}$ is a normal subgroup of $G$. 
Using the
induction hypothesis we obtain 
\begin{multline*}
[G,\lcs H{i-2},H]\leqslant [\lcs G{i-1},H]=[\lcs H{i-1}\lcs
Gi,H]\\ =[\lcs H{i-1},H][\lcs Gi,H]\leqslant
\lcs H{i}\lcs G{i+1}
\end{multline*}
and
\begin{multline*}
[H,G,\lcs H{i-2}]\leqslant[G',\lcs H{i-2}]=[H'\lcs G3,\lcs H{i-2}]\\
=[H',\lcs H{i-2}][\lcs G3,\lcs H{i-2}]\leqslant
\lcs H{i}\lcs G{i+1}.
\end{multline*}
Using the Three Subgroups Lemma (see~\cite{Huppert}~III.1.10(b)), we
obtain 
$$
[\lcs H{i-1},G]=[\lcs H{i-2},H,G]\leqslant\lcs Hi\lcs G{i+1},
$$ 
and hence our statement is correct.

Let us prove that $\lcs Gi=\lcs Hi$ for all $i\geqslant 2$. If
the nilpotency class of $G$ is $c$, that is $\lcs G{c+1}=1$, then
$\lcs G{c+1}=\lcs H{c+1}=1$. If $\lcs
G{i+1}=\lcs H{i+1}$ for some $i$, such that $3\leqslant i+1\leqslant c+1$, then, by
the result of the previous paragraph,
$$
\lcs Gi=\lcs Hi\lcs G{i+1}=\lcs Hi\lcs H{i+1}=\lcs Hi.
$$
Using induction, we obtain $\lcs Gi=\lcs Hi$ for all $i\geqslant 2$.
The normality of $H$ is an easy consequence of the fact that 
$G'=H'\leqslant H$.
\end{proof}

\begin{corollary}\label{2gen3gen}
Let $G$ be a finite $p$-group.
\begin{enumerate}
\item[(i)] If $G'/\lcs{G}{3}$ is cyclic of order $p$, then $G$ has  a
$2$-generator normal subgroup $H$, such that
$\lcs{G}{i}=\lcs{H}{i}$ for all $i\geqslant 2$.
\item[(ii)] If $G'/\lcs{G}{3}$ is elementary abelian of order $p^2$, then
$G$ has  a $3$-generator normal subgroup $H$, such that
$\lcs{G}{i}=\lcs{H}{i}$ for all $i\geqslant 2$. 
\end{enumerate}
\end{corollary}
\begin{proof}
(i)\ Suppose that $G'/\lcs G3=\gen{[a,b]\lcs G3}$ for some $a,\ b\in
G$, and set $H=\gen{a,\ b}$. 
As we have $H'\lcs G3=G'$, Lemma~\ref{gh} implies that
$H$ is a normal subgroup and $\lcs{G}{i}=\lcs{H}i$ for all $i\geqslant
2$. 

(ii)\ Suppose that $G'/\lcs{G}{3}$ is elementary abelian of order
$p^2$, and suppose that $G'/\lcs{G}{3}=\gen{[a,b]\lcs{G}{3},\
[c,d]\lcs{G}{3}}$ for some $a,\ b,\ c,\ d\in G$. Select a subgroup $H$
in $G$ as follows. If $[a,c]$, $[a,d]$, $[b,c]$, $[b,d]$ are all in $\lcs
G3$ then let $H=\gen{a,\ bc,\ d}$. Otherwise suppose without loss of
generality that $[a,c]\equiv [a,b]^\alpha[c,d]^\beta\bmod \lcs{G}{3}$ for some
$\alpha$ and $\beta$, such that 
$0\leqslant \alpha,\ \beta\leqslant p-1$, and at least one of
$\alpha$ and $\beta$ is non-zero. If $\alpha\neq 0$, then set
$H=\gen{a,\ c,\ d}$, otherwise set $H=\gen{a,\ b,\ c}$. It is easy to see
that $H'\lcs G3=G'$, and so, using Lemma~\ref{gh}, we obtain that
$H$ is a normal subgroup and $\lcs{G}{i}=\lcs{H}i$ for all $i\geqslant
2$. 
\end{proof}

\section{Proof of Theorem~\ref{3pgr}}\label{sec2}

In this section we prove Theorem~\ref{3pgr}.

Suppose first that $G$ is a finite $p$-group, such that $|G'/G''|=p^3$
and $G''\neq 1$.
If the quotient $G'/\lcs G3$ is cyclic, then 
Lemma~2.1 of Blackburn~\cite{Blackburn54} implies that
\begin{equation}\label{eq2}
G''=[G',G']=[G',\lcs G3]\leqslant \lcs G5,
\end{equation}
and there is a chain 
\begin{equation}\label{longchain}
G>G'=\lcs G2>\lcs G3>\lcs G4>\lcs G5\geqslant G''>1
\end{equation}
of normal subgroups.  In particular, if $|G'/\lcs G3|=p$,
then~\eqref{eq2} and~\eqref{longchain} imply that $G''=\lcs
G5$; similarly if $|G'/\lcs G3|=p^2$, then~\eqref{eq2}
implies that $G'/\lcs G3$
must be elementary abelian.

Now assume that $|G'/\lcs G3|=p^2$; we show that this can only
happen when $p=2$. 
By Corollary~\ref{2gen3gen}, there is a 3-generator subgroup $H$ of $G$,
such that $\lcs Gi=\lcs Hi$ for all $i\geqslant 2$. After replacing
$G$ by $H$, we may assume without loss of generality that 
$G=\gen{a,\ b,\ c}$ for some $a,\ b,\ c\in
G$. Moreover, from~\eqref{normsgps}  it follows that
$G''=\lcs G4$, and hence we may suppose that $G$ has nilpotency class~4. As
$G'/\lcs G3$ is elementary abelian of order $p^2$, we have that there
are some $\alpha$, $\beta$, and $\gamma$ not all zero, such that $0\leqslant \alpha,\ \beta,\
\gamma\leqslant p-1$ and 
$$
[a,b]^\alpha[a,c]^\beta[b,c]^\gamma\equiv 1\bmod\lcs G3.
$$
If $\alpha=\beta=0$, then $[b,c]^\gamma\equiv 1\bmod\lcs G3$, that is
$[b,c]\in\lcs G3$. If $\alpha=0$ and $\beta\neq 0$, then we obtain
$[a^\beta b^\gamma, c]\equiv 1\bmod\lcs G3$. If we replace
$a^\beta b^\gamma\replace a$, then we obtain that in the new
generating set $[a,c]\in\lcs G3$. Similarly, if $\alpha\neq 0$ and
$\beta=0$, then we replace $a^{-\alpha}c^\gamma\replace a$, and
obtain that after the substitution $[a,b]\in\lcs G3$. If $\alpha\neq 0$, and
$\beta\neq 0$, then we replace
$a^{\beta/\alpha}b^{\gamma/\alpha}\replace a$ and
$bc^{\beta/\alpha}\replace b$.
Then it is easy to see that in the new generating set $[a,b]\in\lcs
G3$. After possibly reordering the generators, we may suppose without
loss of generality that $G$ is generated by three elements $a$, $b$,
and $c$, such that $G'/\lcs G3=\gen{[a,b]\lcs G3,\ [a,c]\lcs G3}$, and
moreover $[b,c]\in\lcs G3$. Note that in this case
$[a,b,c]\equiv[a,c,b]\bmod \lcs G4$ also holds. 
Then 
$G''\neq 1$ and $\lcs G5=1$ imply that
$$
[[a,b],[a,c]]=[a,b,a,c][a,b,c,a]^{-1}\neq 1
$$
and
$$
[[a,c],[a,b]]=[a,c,a,b][a,c,b,a]^{-1}\neq 1.
$$
If $[a,b,a]\in\lcs G4$, then $[a,b,a,c]\in\lcs G5$, and hence
$[a,b,a,c]=1$. Similarly $[a,b,c]\in\lcs G4$, implies that $[a,b,c,a]=1$;
therefore at least one of the elements $[a,b,a]$ and $[a,b,c]$ does
not lie in $\lcs G4$. 
Similarly, at least one of 
$[a,c,a]$ and $[a,b,c]$ must also lie outside $\lcs G4$.

First we assume that $[a,b,c]\in\lcs G4$. 
In this case we must have $[a,c,a]\not\in\lcs G4$ and
$[a,b,a]\not\in\lcs G4$. 
As $\lcs G3/\lcs G4$ is cyclic of order $p$,  there is some
$\alpha$, such that $0\leqslant \alpha\leqslant p-1$ and $[a,bc^\alpha,a]\equiv 1\bmod \lcs
G4$, and we carry out the replacement $bc^\alpha\replace b$. In the
new generating set $[b,c]\in\lcs G3$
still holds, and, in addition, we obtain $[a,b,a]\in\lcs G4$.

So without loss of generality we  assume that $[a,b,a]\in\lcs G4$ and
$[a,b,c]\not\in\lcs G4$. In this case $[a,b,c,a]=[a,c,b,a]\neq 1$, in
other words
$a\not\in \cent{G}{\lcs G3}$.   
On the other hand,     $[a,b,b,a]=[a,b,a,b]$, and hence $[a,b,b,a]=1$. If
$[a,b,b]\not\in\lcs G4$, then $\lcs G3=\gen{[a,b,b],\lcs G4}$, and so
$a\in\cent G{\lcs G3}$, which is impossible; therefore
$[a,b,b]\in\lcs G4$. 
If $[a,c,a]\not\in\lcs G4$ then there is some $\alpha\neq 0$,
such that $[a,c,ab^\alpha]\in\lcs G4$; in this case we let $ab^\alpha\replace a$ and obtain 
$[a,c,a]\in\lcs G4$. 
In the new generating set $[b,c]\in\lcs G3$ and $[a,b,a],\ [a,b,b]\in\lcs
G4$ still hold. Then
\begin{multline*}
1=[[a,b],a,c][[a,c],[a,b]][c,[a,b],a]\\
={[a,b,a,c]}[a,c,a,b][a,b,c,a]^{-1}[a,b,c,a]^{-1}=[a,b,c,a]^{-2},
\end{multline*}
that is, $[a,b,c,a]^2=1$, and hence $p=2$. 
This completes the proof of Theorem~\ref{3pgr}.

\section{The proof of Theorem~\ref{decomp2}}\label{sec3}

In the previous section we proved Theorem~\ref{3pgr}, and hence
we know that in a group $G$ the conditions of Theorem~\ref{decomp2}
imply that $|G'/\lcs{G}{3}|=p$.  Thus, according to Corollary~\ref{2gen3gen}, $G$ has a 2-generator subgroup $H$, such
that for all $i\geqslant 2$ we have $\lcs Gi=\lcs Hi$. We use
this subgroup to obtain the desired factorisation. First we show that
we can choose a generating set which satisfies some extra conditions.

\begin{lemma}\label{2gen}
Let $G$ be a $2$-generator, finite $p$-group, such that 
$|G'/G''|=p^3$, $|G'/\lcs{G}{3}|=p$, and $G''\neq 1$. 
Then  generators $a$ and $b$ of $G$
can be chosen, such that the following hold:
\begin{enumerate}
\item[(i)] $\lcs{G}{2}/\lcs{G}{3}=\gen{[b,a]\lcs{G}{3}}$;
\item[(ii)] $\lcs{G}{3}/\lcs{G}{4}=\gen{[b,a,a]\lcs{G}{4}}$ and $[b,a,b]\in\lcs{G}{4}$;
\item[(iii)] $\lcs{G}{4}/\lcs{G}{5}=\gen{[b,a,a,a]\lcs{G}{5}}$ and $[b,a,a,b]\in\lcs{G}{5}$;
\item[(iv)] $\lcs{G}{5}/\lcs{G}{6}=\gen{[b,a,a,a,b]\lcs{G}{6}}$ and $[b,a,a,a,a]\in\lcs{G}{6}$. 
\end{enumerate}
\end{lemma}
\begin{proof}
We may suppose without loss of generality that $G$ has class~5. 
As noticed in the introduction, 
our conditions imply that the factors
$G'/\lcs G3$, $\lcs G3/\lcs G4$, and $\lcs G4/\lcs G5$
are cyclic with order $p$. Using the argument presented by
Blackburn~\cite{Blackburn54} in Lemma~2.9, we can choose the generating
set $\{a,\ b\}$, so that properties~(i)-(iii) hold. It follows
from~\eqref{eq2} and~\eqref{longchain} that $G''=\lcs G5$, and
$G''=\gen{[[b,a,a],[b,a]]}$. As the element $[[b,a,a],[b,a]]$ is
central and has order $p$, we have $|\lcs G5|=p$, and using
Blackburn's argument on page~89, the set $\{a,\ b\}$ can be chosen so
that the additional property~(iv) also holds.
\end{proof}

\begin{lemma}\label{decomp}
Let $p\geqslant 3$ and $G$ be a finite $p$-group, such that $|G'/G''|=p^3$ and
$G''\neq 1$. Then $G$ has a minimal generating set
$\{a,b,u_1,u_2,\ldots,u_r\}$, such that
\begin{enumerate}
\item[(i)] $H=\bigl<a,b\bigr>$ is a normal subgroup of $G$, such that
$\lcs Hi=\lcs Gi$ for all $i\geqslant 2$; further, $a$ and $b$ are as in Lemma~{\rm \ref{2gen}};
\item[(ii)] $[a,u_i]\in\lcs{G}{5}$ for all $u_i$;
\item[(iii)] $[b,u_i]\in\lcs{G}{4}$ for all $u_i$;
\item[(iv)] $[u_i,u_j]\in\lcs{G}{5}$ for all $u_i$ and $u_j$.
\end{enumerate}
In particular, $u_1,\ldots,u_r\in\cent G{G'}$.
\end{lemma}
\begin{proof}
First recall Hall's theorem that $|G''|=p$, and
 so~\eqref{longchain} implies that $G$ has class 5.
Select $a,\ b\in G$, such that the subgroup $H=\gen{a,b}$ 
 and its generators are as in Lemma~\ref{2gen}. 
It is easy to see that $a,\ b$ are linearly 
independent modulo the Frattini subgroup of $G$. Therefore they can be viewed 
as elements of a
minimal generating set $\{a,b,u_1,\ldots,u_r\}$. Now suppose that
 for each $i\in\{1,\ldots,r\}$, 
$[u_i,a]\equiv[b,a]^{\alpha_i}$ and $[u_i,b]\equiv[b,a]^{\beta_i}$ modulo
 $\lcs{G}{3}$ with some $\alpha_i,\ \beta_i\in\{0,\ldots,p-1\}$. Then $[u_ib^{-\alpha_i} a^{\beta_i},b]\in\lcs{G}{3}$ and also 
$[u_ib^{-\alpha_i} a^{\beta_i},a]\in\lcs{G}{3}$. If we perform the
 replacement  $u_ib^{-\alpha_i} a^{\beta_i}\replace u_i$, then
it is easy to see that $\{a,b,u_1,\ldots,u_r\}$ is also a minimal 
generating set for $G$ and 
$$
[a,u_i],\ [b,u_i]\in\lcs{G}{3}\quad\mbox{for all}\quad i\in\{1,\ldots,r\}.
$$

Now suppose that for all $i\in\{1,\ldots,r\}$ we have
$$
[u_i,a]\equiv[b,a,a]^{\alpha_i}[b,a,a,a]^{\beta_i}\bmod \lcs{G}{5}
$$ 
for some $\alpha_i,\ \beta_i\in\{0,\ldots,p-1\}$. Then computing modulo
$\lcs{G}{5}$ we obtain
\begin{multline*}
[u_i[b,a]^{-\alpha_i}[b,a,a]^{-\beta_i},a]\\
={[u_i[b,a]}^{-\alpha_i},a][u_i[b,a]^{-\alpha_i},a,[b,a,a]^{-\beta_i}][[b,a,a]^{-\beta_i},a]]\\
\equiv{[u_i,a]}[u_i,a,[b,a]^{-\alpha_i}][[b,a]^{-\alpha_i},a][[b,a,a]^{-\beta_i},a]\\\equiv{[u_i,a]}[b,a,a]^{-\alpha_i}[b,a,a,a]^{-\beta_i}\equiv
1.
\end{multline*}
If we replace $u_i[b,a]^{-\alpha_i}[b,a,a]^{-\beta_i}\replace u_i$,
then $[u_i,a]\in\lcs{G}{5}$. Since the images of the $u_i$ over the Frattini subgroup did not change, 
the set $\{a,b,u_1,\ldots,u_r\}$ is still a minimal generating system
for $G$. We show that this generating set satisfies the
properties required by the lemma. 

We claim that  $[u_i,b]\in\lcs{G}{4}$ for all 
$i\in\{1,\ldots,r\}$. To prove this we observe that
$$
1=[b,a,u_i[u_i,b]][u_i,b,a[a,u_i]][a,u_i,b[b,a]]=[b,a,[u_i,b]][b,a,u_i][u_i,b,a],
$$
and thus 
\begin{multline}\label{eq}
[b,a,u_i]=
[u_i,b,a]^{-1}[b,a,[u_i,b]]^{-1}=[[b,u_i]^{-1},a]^{-1}[[u_i,b],[b,a]]\\
=[b,u_i,a][[u_i,b],[b,a]].
\end{multline}
In particular, $[b,a,u_i]\in\lcs{G}{4}$. Now consider
$$
1=[[b,a],a,u_i[u_i,[b,a]]][u_i,[b,a],a[a,u_i]][a,u_i,[b,a][b,a,a]]=[b,a,a,u_i].
$$
We can obtain similarly $[b,a,a,a,u_i]=1$.
As $[u_i,b]\in\lcs{G}{3}$, $[u_i,b]\equiv[b,a,a]^{\varepsilon_i}$
modulo $\lcs{G}{4}$ for some 
$\varepsilon_i\in\{0,\ldots,p-1\}$.
The Hall-Witt identity implies that
\begin{multline*}
1=[[b,a],u_i,b[b,[b,a]]][b,[b,a],u_i[u_i,b]][u_i,b,[b,a][b,a,u_i]]\\
=[[b,a],u_i,b][[u_i,b],[b,a]].
\end{multline*}
Using (\ref{eq}) we get
$[b,a,u_i,b]=[b,a,a,a,b]^{-\varepsilon_i}$. Moreover,  
$$
[[u_i,b],[b,a]]=[[b,a,a]^{\varepsilon_i},[b,a]]=[b,a,a,a,b]^{-\varepsilon_i},
$$
and thus $[b,a,a,a,b]^{-2\varepsilon_i}=1$, from which it follows that
$\varepsilon_i=0$, in other words $[u_i,b]\in\lcs{G}{4}$. 

We now prove that $u_1,\ldots,u_r\in \cent G{G'}$. We have already
seen that $[b,a,a]$, $[b,a,a,a]$ are centralised by the $u_i$, so it
suffices to prove that $[b,a,u_i]=1$ for all $i\in\{1,\ldots,r\}$. This
is clear because
$$
1=[b,a,u_i[u_i,b]][u_i,b,a[a,u_i]][a,u_i,b[b,a]]=[b,a,u_i].
$$

It remains to show that $[u_i,u_j]$ lies 
in $\lcs{G}{5}$ for all $i,\ j\in\{1,\ldots,r\}$. It easily follows using
the Hall-Witt identity that 
$[u_i,u_j,a]=1$ and $[u_i,u_j,b]=1$, therefore $[u_i,u_j]\in Z(H)\cap
H'=\lcs{H}{5}=\lcs{G}{5}$. The proof is complete.
\end{proof}

\noindent{\sc Proof of Theorem~\ref{decomp2}.}
Choose a generating set $\{a,b,u_1,u_2\ldots,u_r\}$ 
for $G$ as in the previous lemma. In the first stage of the
proof we show that this generating set can be modified so that, in
addition to the properties required by Lemma~\ref{decomp}, one
of the following holds:
\begin{enumerate}
\item[(a)] $u_1,\ldots,u_r\in\cent G a$; or
\item[(b)] $u_2,\ldots,u_r\in\cent G{\gen{a,u_1}}$.
\end{enumerate}

If
$u_1,\ldots,u_r\in \cent Ga$ then~(a) holds and we are done. Suppose that there is at least one $u_i$ which
does not centralise $a$. Without loss of generality we may assume that
$[u_1,a]=[b,a,a,a,b]$. If $[u_i,a]=[b,a,a,a,b]^{\alpha_i}$ for some
$i\in\{2,\ldots,r\}$,  then let $u_iu_1^{-\alpha_i}\replace
u_i$. In this way
we obtain a generating set $\{a,b,u_1,\ldots,u_r\}$,  such that
$[u_1,a]=[b,a,a,a,b]$ and $\gen{u_2,\ldots,u_r}\leqslant \cent Ga$.

If $u_2,\ldots,u_r$ centralise $u_1$, then~(b) holds and we are done. 
We assume without loss of generality that
$[u_2,u_1]=[b,a,a,a,b]$. If $[u_i,u_1]=[b,a,a,a,b]^{\beta_i}$ for some
$i\in\{3,\ldots,r\}$, then
let $u_iu_2^{-\beta_i}\replace u_i$. In this way we obtain a generating
set, such that $u_2,\ldots,u_r$ centralise
$a$, and $u_3,\ldots,u_r$ centralise $u_1$. Repeating this process, we
construct a generating set $\{a,b,u_1,\ldots,u_k,\ldots,u_r\}$, such that 
\begin{enumerate}
\item $[u_1,a]=[b,a,a,a,b]$; 
\item $u_2,\ldots,u_r$ centralise $a$; 
\item $[u_{i+1},u_i]=[b,a,a,a,b]$ for all $i\in\{1,\ldots,k-1\}$;
\item $[u_{k+1},u_k]=1$;
\item $u_{i+2},\ldots,u_r$ centralise $u_i$ for all $i\in\{1,\ldots,k\}$.
\end{enumerate}
Now if $k$ is even then substitute $au_2u_4\cdots u_k\replace a$. After this change property~(a) holds. 
If $k$ is odd then replace $u_1u_3\cdots u_k\replace u_1$;
in this case property~(b) holds.

We continue with the second stage of the proof. Suppose that
the generating set $\{a,b,u_1,\ldots,u_r\}$ is as in Lemma~\ref{decomp} and, in addition,
property~(a) holds.
First assume that all the $u_i$ centralise $b$ modulo $\lcs{G}{5}$. If
$[u_i,b]=[b,a,a,a,b]^{\gamma_i}$ for some $i\in\{1,\ldots,r\}$ and $\gamma_i\in\{0,\ldots,p-1\}$, then let
$u_i[b,a,a,a]^{-\gamma_i}\replace u_i$. Then $H=\gen{a,b}$ and
$U=\gen{u_1,\ldots,u_r}$ satisfy the assertions of the theorem.

Suppose that some of the $u_i$ do not
centralise $b$ modulo $\lcs{G}{5}$, and assume without loss of generality that 
$[u_1,b]=[b,a,a,a][b,a,a,a,b]^{\gamma_1}$. Perform the substitution  
$u_1[b,a,a,a]^{-\gamma_1}\replace u_1$ to obtain $[u_1,b]=[b,a,a,a]$. If
$[u_i,b]\equiv[b,a,a,a]^{\gamma_i}\bmod \lcs{G}{5}$ with some
$i\in\{2,\ldots,r\}$, then substitute 
$u_iu_1^{-\gamma_i}\replace u_i$. After this there is some $\delta_i$,
such that $0\leqslant \delta_i\leqslant
p-1$ and $[u_i,b]=[b,a,a,a,b]^{\delta_i}$; then
replace $u_i[b,a,a,a]^{-\delta_i}\replace u_i$. This way we obtain 
$[u_1,b]=[b,a,a,a]$ and, moreover,  
$\gen{u_2,\ldots,u_r}\leqslant \cent G{b}$. If $u_2,\ldots,u_r$ centralise
$u_1$, then choose $H=\gen{a,b,u_1}$ and $U=\gen{u_2,u_3,\ldots,u_r}$
and we are done. Suppose that this is not the case and
$[u_2,u_1]=[b,a,a,a,b]$. Then, as in the first part of the proof,
select a generating set $\{a,b,u_1,\ldots,u_k,\ldots,u_r\}$,
such that the following additional properties hold:
\begin{enumerate}
\item $[u_1,b]=[b,a,a,a]$; 
\item $u_2,\ldots,u_r$ centralise $b$; 
\item $[u_{i+1},u_i]=[b,a,a,a,b]$ for all $i\in\{1,\ldots,k-1\}$;
\item $[u_{k+1},u_k]=1$;
\item $u_{i+2},\ldots,u_r$ centralise $u_i$ for all $i\in\{1,\ldots,k\}$.
\end{enumerate}
If $k$ is even then set 
$$
H=\gen{a,b,u_1u_3\cdots u_{k-1},u_2u_4\cdots
u_k}
$$
and 
$$
U=\gen{u_2,u_3,\ldots,u_{k-1},u_{k+1},\ldots, u_r}.
$$
If $k$ is odd then let $H=\gen{a,b,u_1u_3\cdots u_{k}}$  and
$U=\gen{u_2,u_3,\ldots,u_r}$. In both cases the subgroups $H$ and $U$
are as required.

In the case of property~(b), we consider the group
$G_1=\gen{a,b,u_2,\ldots,u_r}$ and choose subgroups $H_1$ and $U_1$
according to the process described in the previous paragraph. Then
note that $H_1$ and $U_1$ satisfies the prescribed conditions.
Moreover $H_1$ can be generated by at most four elements. For $G$ we
can choose the subgroups $H=\gen{H_1,u_1}$ and $U=U_1$.
\hfill$\Box$

\bigskip

The following example shows that the number ``5'' in
Theorem~\ref{decomp2} is the best possible. This construction can be
generalised, and it is not difficult to see that similar examples exist
for all $p$. 

\begin{example}
Consider the pro-5-group $G$ given by the pro-5-presentation
\begin{multline*}
\{a,\,b,\,u_1,\,u_2,\,u_3\,|\,a^5,\,b^5,\,u_1^5,\,u_2^5,\,u_3^5,\\
{[b,a,b]},\,[b,a,a,a,a],\,
[b,a,a,a,b][a,u_1],
[a,u_2],[a,u_3],
\\{[b,u_1]},\,[b,a,a,a][b,u_2],\,[b,u_3],\,
[u_1,u_2],\,[u_1,u_3],\,[b,a,a,a,b][u_2,u_3]\}.
\end{multline*}
Then, using the ANU~$p$-Quotient Program~\cite{pq,pqart}, it is easy to
see that $G$ is a finite 5-group and $\lcs{G}{5}=G''\neq 1$. 
Suppose that $G=HU$ is a
factorisation of $G$ as in the theorem. Then $U$ centralises $H$, and
in particular, $U\leqslant \cent{G}{G'}$. Using a
computer algebra system, such as {\sf GAP}~\cite{gap} or {\sc Magma}~\cite{magma}, it is easy
to compute that $\cent G{G'}=\gen{u_1,u_2,u_3,[b,a,a,a]}$, and that no
subgroup of $G$ generated by less than 5 generators
can be taken for $H$ in Theorem~\ref{decomp2}.
\end{example}

In Theorem~\ref{decomp2} the subgroup $U$ satisfies
$|U'|\leqslant p$. The non-abelian $p$-groups with this property were classified by
S.~R.~Blackburn~\cite{sblackburn}. Unfortunately, the isomorphism
types of $H$ and $U$ are not uniquely 
determined by the isomorphism type of $G$.
The following example illustrates this fact. 

\begin{example}
Let $p\geqslant 5$ and let $G$ denote the pro-$p$-group given by the
pro-$p$-presentation
\begin{multline*}
\{a,\,b,\,u_1,\,u_2,\,u_3 |\,
a^p,\,b^p,\,u_1^{p^3},\,u_2^{p^2},\,u_3^{p^2},\,
[b,a,b],\,[b,a,a,a,a],\\
[b,a,a,a,b][a,u_1],
[a,u_2],\,
[b,u_1],\,[b,u_2],\,[u_1,u_2],\\
[u_3,a],\,[u_3,b],\,[u_3,u_1],\,[b,a,a,a,b][u_3,u_2]\}.
\end{multline*}
Then $G$ has the obvious factorisation $G=H_1U_1$, where
$H_1=\gen{a,b,u_1}$ and $U_1=\gen{u_2,u_3}$. The group $G$ also admits
a factorisation $G=H_2U_2$, where $H_2=\gen{au_3,b,u_1}$ and
$U_2=\gen{u_1u_2^{-1},u_3}$. It is easy to see that $H_1\not\cong H_2$
and $U_1\not\cong U_2$. 
\end{example}

\section{A characterisation of the derived subgroup}\label{sec4}

The following lemma was already known to Burnside. Its proof
is an easy exercise, and can also be found in Huppert~\cite{Huppert}~III.7.8.

\begin{lemma}\label{huplemma}
In a finite $p$-group $G$, if $Z(G')$ is cyclic then so is $G'$. 
\end{lemma}

Suppose that $G$ is a $p$-group for some odd $p$, such that
$|G'/G''|=p^3$ and $G''\neq 1$.  As $|G''|=p$, the subgroup $G'$ has
order $p^4$ and its derived subgroup
$G''$ is cyclic
with order $p$. By the previous lemma $Z(G')$ cannot be
cyclic. The following result gives more information on the
structure of $G'$.

\begin{lemma}
The quotient $G'/G''$ is elementary abelian.
\end{lemma}
\begin{proof}
Recall that Hall's theorem
implies that $\lcs{G}{6}=1$. 
Using Corollary~\ref{2gen3gen}, assume that $G$ is generated
by two elements 
$a$ and $b$ which are chosen as in Lemma~\ref{2gen}. Then $G'/G''$ is generated by the images of $[b,a]$, $[b,a,a]$ and
$[b,a,a,a]$. Since the centre of $G'$ is
$\gen{[b,a,a,a],G''}$, we must have
$[b,a,a,a]^p=1$ by Lemma~\ref{huplemma}. 

Suppose that
$[b,a]^p\not\equiv 1\bmod\lcs{G}{4}$. Then 
$[b,a]^p\lcs{G}{4}$ generates the factor $\lcs{G}{3}/\lcs{G}{4}$ and in
particular $[[b,a]^p,[b,a]]\neq 1$, which is clearly impossible. 
Suppose now that $[b,a]^p\not\equiv 1\bmod\lcs{G}{5}$. Then
$$
[[b,a]^p,b]\equiv[b,a,b]^p=1\bmod (N')^p\lcs{N}{p},
$$
where
$N=\gen{[b,a],[b,a,b]}$.
This yields $[[b,a]^p,b]=1$, which is a contradiction. Now suppose
that $[b,a,a]^p\not\equiv 1\bmod\lcs{G}{5}$. Then 
$$
[[b,a,a]^p,b]\equiv[b,a,a,b]^p=1\bmod(N')^p\lcs{N}{p},
$$
where
$N=\gen{[b,a,a],[b,a,a,b]}$.
Again, this leads to a contradiction.
\end{proof}

Our last main result is a characterisation of $G'$. For odd primes let $X_{p^3}$ and
$Y_{p^3}$ denote the non-abelian $p$-groups of order $p^3$ and exponent
$p$ and $p^2$, respectively. The symbol $C_p$ denotes the cyclic group
of order $p$.

\begin{theorem}\label{isomcorol}
If $p\geqslant 3$  and $G$ is a finite $p$-group, such that
$|G'/G''|=p^3$ and $G''\neq 1$, then $G'$ is isomorphic to $X_{p^3}\times C_p$ or to $Y_{p^3}\times
C_p$.
\end{theorem}
\begin{proof}
Recall that by Hall's theorem $|G'|=p^4$. 
For $p\geqslant 5$ the list of groups with order $p^4$ can be found in Huppert~\cite{Huppert}~III.12.6. For $p=3$ one can find this list as part of {\sf
  GAP}~\cite{gap} or {\sc Magma}~\cite{magma}. 
It is easy to see that the only groups which satisfy the
conditions on $G'$ are $X_{p^3}\times C_p$ and  $Y_{p^3}\times
C_p$.
\end{proof}

\begin{example}
Let $G$ be a group of maximal class of order $p^6$ for $p\geqslant 5$ with
degree of commutativity $0$. Then $|G'/G''|=p^3$ and 
by Theorem 3.2 of Blackburn~\cite{Blackburn54} $G'\cong X_{p^3}\times
C_p$. An example for such a group is the 
pro-$p$-group described by the pro-$p$-presentation
$$
G=\{a,b\mid a^p,b^p,[b,a,b],[b,a,a,a,a]\}.
$$

If $p=3$ then the pro-3-group described by the pro-3-presentation
$$
\{a,b\mid a^9,b^9,[a,b]^3,[b,a,b],[b,a,a,a,a]\}
$$
contains $X_{27}\times C_3$ as derived subgroup. This can easily be checked
using the $p$-Quotient Program (\cite{pq,pqart}).
\end{example}

\begin{example}
If $p\geqslant 3$ and  $G$ denotes the pro-$p$-group given by the
pro-$p$-presentation 
$$
\{a,b\mid a^{p^2},b^{p^2},[b,a]^p=[b,a,a,a,b],[b,a,b],[b,a,a,a,a]\},
$$
then $G'\cong Y_{p^3}\times C_p$.
\end{example}

\begin{corollary}\label{emaxclass}
If $p\geqslant 5$ and $G$ is a finite $p$-group, such that
$G'\cong X_{p^3}\times C_p$, then $G^p\leqslant Z(G)$. If $p\geqslant 3$ and $G$ is a finite $p$-group, such that
$G'\cong Y_{p^3}\times C_p$, then $G^{p^2}\leqslant Z(G)$.
\end{corollary}
\begin{proof}
We only prove the first statement; the proof of the second is very
similar.
It is enough to prove that $u^p\in Z(G)$ for all $u\in G$.
So let $u\in G$ and notice that $[v,u]\in G'$ for all $v\in G$. 
By the
collection formula
$$
[v,u^p]\equiv[v,u]^p=1\bmod (N')^p\lcs{N}{p}\quad\mbox{where}\quad N=\gen{u,[u,v]}.$$
If $p\geqslant 5$ then $(N')^p\lcs{N}{p}=1$ therefore
$[v,u^p]=1$.
\end{proof}

\section{Acknowledgement}
Much of the research presented in this paper was carried out while I
was a PhD student at The Australian National University in
Canberra. I am grateful to my supervisor, Mike Newman, for his
many suggestions; and to Cheryl Praeger and Eamonn O'Brien for their
helpful comments on a draft.

\bibliographystyle{plain}

\end{document}